\documentclass[graybox]{svmult}

\usepackage{mathptmx}       
\usepackage{helvet}         
\usepackage{courier}        
\usepackage{type1cm}        
\usepackage{makeidx}         
\usepackage{graphicx}        
\usepackage{multicol}        
\usepackage[bottom]{footmisc}

\usepackage{epstopdf}
\usepackage{overpic}

\makeindex             
\usepackage{amssymb,amsmath}
\usepackage[chapter]{algorithm}
\usepackage{algpseudocode}
\usepackage{algorithmicx}

\newcommand{\bp}{\mathbf{p}}

\newcommand{\bu}{\mathbf{u}}
\newcommand{\bv}{\mathbf{v}}

\newcommand{\bx}{\mathbf{x}}

\newcommand{\by}{\mathbf{y}}
\newcommand{\bz}{\mathbf{z}}
\newcommand{\bA}{\mathbf{A}}
\newcommand{\bB}{\mathbf{B}}
\newcommand{\bC}{\mathbf{C}}

\newcommand{\bF}{\mathbf{F}}
\newcommand{\bG}{\mathbf{G}}
\newcommand{\bH}{\mathbf{H}}
\newcommand{\bK}{\mathbf{K}}

\newcommand{\bP}{\mathbf{P}}

\newcommand{\bS}{\mathbf{S}}

\newcommand{\bLambda}{\boldsymbol{\Lambda}}

\graphicspath{{FIGURES/}}

\begin{document}

\title*{Data-driven approximations of dynamical systems operators for control}
\author{Eurika Kaiser, J.~Nathan Kutz, and Steven L. Brunton}
\institute{Eurika Kaiser \at Mechanical Engineering, University of Washington, Seattle, WA, 98195, \email{eurika@uw.edu}
\and J.~Nathan Kutz \at Applied Mathematics, University of Washington, Seattle,
WA, 98195 \email{kutz@uw.edu}
\and Steven~L. Brunton \at Mechanical Engineering, University of Washington, Seattle, WA, 98195, \email{sbrunton@uw.edu}}
%
%
\maketitle

\abstract*{
The Koopman and Perron Frobenius transport operators are fundamentally changing how we approach dynamical systems, providing linear representations for even strongly nonlinear dynamics.  
Although there is tremendous potential benefit of such a linear representation for estimation and control, transport operators are infinite-dimensional, making them difficult to work with numerically.  
Obtaining low-dimensional matrix approximations of these operators is paramount for applications, and the dynamic mode decomposition has quickly become a standard numerical algorithm to approximate the Koopman operator.  
Related methods have seen rapid development, due to a combination of an increasing abundance of data and the extensibility of DMD because of its simple framing in terms of linear algebra.  
In this chapter, we review key innovations in the data-driven characterization of transport operators, providing a high-level and unified perspective.  
We emphasize important recent developments around sparsity and control, in particular, and discuss emerging methods in big data and machine learning.  	
}

\abstract{
The Koopman and Perron Frobenius transport operators are fundamentally changing how we approach dynamical systems, providing linear representations for even strongly nonlinear dynamics.  
Although there is tremendous potential benefit of such a linear representation for estimation and control, transport operators are infinite-dimensional, making them difficult to work with numerically.  
Obtaining low-dimensional matrix approximations of these operators is paramount for applications, and the dynamic mode decomposition has quickly become a standard numerical algorithm to approximate the Koopman operator.  
Related methods have seen rapid development, due to a combination of an increasing abundance of data and the extensibility of DMD based on its simple framing in terms of linear algebra.  
In this chapter, we review key innovations in the data-driven characterization of transport operators for control, providing a high-level and unified perspective.
%
We emphasize important recent developments around sparsity and control, and discuss emerging methods in big data and machine learning.  
}


%


%
%
\section{Introduction}
\label{Sec:Introduction}

Data-driven modeling using linear operators  has the potential to transform the estimation and control of strongly nonlinear systems.
Linear operators, such as Koopman and Perron-Frobenius operators, provide a principled linear embedding of nonlinear dynamics, extending the application of standard linear methods to nonlinear systems, and thus, significantly simplifying the control design and reducing the computational burden.
More broadly, data-driven discovery of dynamical systems is undergoing rapid development, driven by the lack of simple, or often entirely unknown, equations and the increasing abundance of high-fidelity data.
There have been recent successes in the discovery of functional representations of nonlinear dynamical systems, e.g.\ using evolutionary optimization techniques~\cite{Bongard2007pnas,Schmidt2009science} and sparse optimization~\cite{Brunton2016pnas}.
However, control based on nonlinear equations is particularly challenging, becoming infeasible for higher-dimensional problems, lacking guarantees, and often requiring problem-tailored formulations.
In contrast, the emerging field of linear operators in dynamical systems seeks to embed nonlinear dynamics in a globally linear representation, providing a compelling mathematical framework for the linear estimation, prediction, and control of strongly nonlinear systems.
The rise of advanced data science and machine learning algorithms, vastly expanded computational resources, and advanced sensor technologies make this a fertile ground for the rapid development of data-driven approximations to these linear operators for control.
In this chapter, we review key innovations, discuss major challenges and promising future directions, and provide a high-level and unified perspective on data-driven approximations of transfer operators for control.

Data-driven system identification has reached a high degree of maturity. 
There exist a plethora of techniques that identify linear and nonlinear systems~\cite{Nelles2013book} based on data, including state-space modeling via the eigensystem realization algorithm (ERA)~\cite{ERA:1985} and other subspace identification methods, Volterra
series~\cite{Brockett1976automatica,maner1994nonlinear}, linear and nonlinear autoregressive models~\cite{Akaike1969annals} (e.g., ARX, ARMA, NARX, and NARMAX), and neural network models~\cite{lippmann1987introduction,draeger1995model,wang2016combined}, to name only a few. 
The reader is referred to~\cite{ljung2010arc,ljung2010ccc} for an compressed overview of identification techniques for linear and nonlinear systems. 
In the machine learning community, manifold learning, e.g. locally linear embedding and self-organizing maps, and non-parametric modeling, e.g.\ Gaussian processes, have been proven to be useful for identifying nonlinear systems~\cite{principe1998ieee,ko2007ieee,kocijan2004acc}. 
Most of these models are considered data-driven~\cite{sjoberg1995automatica}, as they do not impose a specific model structure based on the governing equations. 
However, there is an increasing shift from black-box modeling to inferring unknown physics and constraining models with known prior information.
For instance, the recent sparse identification of nonlinear dynamics (SINDy)~\cite{Brunton2016pnas}, which has been extended to incorporate the effect of control~\cite{Brunton2016nolcos, kaiser2017arxiv_b}, is able to take into account known expert knowledge such as symmetries and conservation laws~\cite{Loiseau2018jfm}. 
Learning accurate nonlinear models is particularly challenging as small deviations in the parameters may produce fundamentally different system behavior.
Importantly, it is possible to control many nonlinear systems using linear models. Examples include weakly nonlinear systems for which models are obtained based on a local linearization of the nonlinear dynamics at a reference point. 
However, the pronounced nonlinearities present in many applications generally require nonlinear control design.
The lack of robustness and stability guarantees except for special cases and the increased computational burden during the on-line phase, which becomes prohibitive for high-dimensional systems, restricts the application of nonlinear control to low-dimensional systems and requires specialized design techniques.
Encoding nonlinear dynamics in linear models through operator-theoretic approaches provides a new opportunity for the control of previously intractable systems. 

Linear embedding theory for nonlinear systems goes back to seminal works by B.O. Koopman in 1931~\cite{Koopman1931pnas} and T. Carleman in 1932~\cite{carleman1932am}. 
Koopman showed that Hamiltonian dynamics can be described through an infinite-dimensional linear operator acting on the Hilbert space of all possible observables that can be measured from the underlying state. 
Closely related, Carleman demonstrated that systems of ordinary differential equations with polynomial nonlinearities can be represented as an infinite-dimensional system of linear differential equations.
Since the introduction of the so-called Carleman linearization, it has been applied to a wide range of problems~\cite{bellman1963qam,brockett1976ieee,kowalski1991nonlinear}
including for the Lyapunov exponent calculation~\cite{andrade1982jmp} and for finding first integrals~\cite{kus1983jpa}.
The emerging Koopman operator perspective provides an alternative direction and generalizes beyond polynomial systems.
The ultimate goal is to learn a globally linear embedding of the nonlinear dynamics such that powerful linear methods become immediately applicable and are useful in a larger domain. 
Of particular interest are the spectral properties of these operators that encode global information and can be related to geometrical properties of the underlying dynamical system. 
The potential to discover global properties of dynamical systems through operator-theoretic methods for diagnostic purposes, i.e. improving understanding of the underlying dynamics, has driven continued efforts to develop improved algorithms.
The confluence of big data, advances in machine learning, and new sensor technologies has further fueled the progress in data-driven methods, facilitating equation-free approximations of these operators. 
These efforts will be discussed below particularly in the context of system identification for control.

We introduce the Koopman and Perron-Frobenius operators by considering the following autonomous nonlinear dynamical system
\begin{equation}\label{Eqn:Dynamics}
\frac{d}{dt} \bx(t) = \bF(\bx(t)),\quad \bx(0) = \bx_0
\end{equation}
with $\bx\in\mathbb{X}\subset\mathbb{R}^n$, initial condition $\bx_0\in\mathbb{X}$, and the flow is denoted by $\bS^t$ so that $\bx(t)=\bS^t(\bx_0)$.
The nonlinear system~\eqref{Eqn:Dynamics} can be equivalently described by infinite-dimensional, linear operators acting on observable or density functions.
The linearity of these operators is appealing; however, their infinite dimensionality poses issues for representation and computation and current research aims to approximate the evolution instead on a finite-dimensional subspace facilitating a finite-dimensional matrix representation~\cite{Brunton2016plosone}.


Let $(\mathbb{X}, \mathfrak{B},\mu)$ be a measure space with state space $\mathbb{X}$, $\sigma$-algebra $\mathfrak{B}$, and measure $\mu$. 
The Koopman operator~\cite{Koopman1931pnas,Koopman1932pnas,Mezic2005nd,Mezic2013arfm} is an infinite-dimensional linear operator that advances measurement functions $f\in L^{\infty}(\mathbb{X})$:
\begin{equation}\label{Eqn:KoopmanOperator}
f(t,\bx_0) = U^t f(\bx_0)= f(\bS^t(\bx_0)). 
\end{equation}
Any set of eigenfunctions of the Koopman operator span an invariant subspace, where the dynamics evolve linearly along these basis directions. 
Thus, the spectral analysis of the Koopman operator is of particular interest and spectral properties have been shown to be related to intrinsic time scales, geometrical properties, and long-term behavior of the underlying dynamical system~\cite{Mezic2017book}.
Koopman eigenfunctions $\phi(\bx)$ associated with a particular eigenvalue $\lambda$ are special or intrinsic observables, which evolve linearly according to
\begin{equation}\label{Eqn:KoopmanEigenfunctions}
U^t \phi(\bx_0)= \phi(\bS^t(\bx_0)) = e^{\lambda t} \phi(\bx_0), 
\end{equation}
and which can generally be discontinuous~\cite{mezic1999chaos,mezic2003cdc,Mezic2004physicad}.
It can be shown~\cite{lasota2013book} for continuously differentiable functions $f$ with compact support, the function ${f}(t,\bx) := U^t[f](\bx)$ satisfies the first-order partial differential equation (PDE):
\begin{equation}
\frac{\partial}{\partial t}{f}(t,\bx) = {\bf F}(\bx)\cdot\nabla {f}(t,\bx) = L_U {f}(\bx),\quad {f}_0:={f}(0,\bx),
\end{equation}
where $L_U$ is the infinitesimal generator of the semigroup of Koopman operators $\{U^t\}_{t\geq 0}$, for which an exponential representation $U^t = e^{L_U t}$ exists.
Smooth Koopman eigenfunctions can then be interpreted as the eigenfunctions of the generator and satisfy
\begin{equation}
\frac{d}{dt}\phi(\bx) = {\bf F}(\bx)\cdot\nabla \phi(\bx)= \lambda \phi(\bx).
\end{equation}
As the Koopman operator evolves measurement functions, this perspective is particularly amenable to data-driven approaches.

The Koopman operator is dual to the Perron-Frobenius operator~\cite{nicolis1995book,lasota2013book,ChaosBook2012,bollt2013book}, i.e. $\langle P^t\rho,f\rangle = \langle \rho,U^t f\rangle$ for any $\rho\in L^1$ and $f\in L^{\infty}$, and as a consequence these share certain properties. For instance, for measure-preserving systems these operators are unitary and exhibit the same point spectrum~\cite{Mezic2004physicad}.
The Perron-Frobenius operator propagates densities $\rho\in L^1(\mathbb{X})$ and is defined as
\begin{equation}
	\int_{B}\, P^t\rho(\bx) \mu(d\bx) =  \int_{\bS^{-t}(B)}\rho(\bx)\,\mu(d\bx)\quad\forall B\in\mathfrak{B}.
\end{equation}
Further, 
\begin{equation}
\rho(t,\bx) = P^t\rho(\bx)= \int_{\mathfrak{B}} \delta(\bx-\bS^t(\bx_0))\rho_0(\bx_0) d\bx_0, 
\end{equation}
where $\delta(\bx-\bS^t(\bx_0))$ represents the deterministic kernel. 
For an invertible system, this becomes $P^t\rho(\bx)= J^{-t}(\bx) \rho(\bS^{-t}(\bx))$,
where $J^{-t}(\bx):=\det(d\bS^{-t}(\bx)/d\bx)$ is the determinant of the Jacobian of $\bS^{-t}(\bx)$;
thus, the density varies inversely with the infinitesimal volume occupied by the trajectories.
For invertible and conservative systems we have $P^t\rho(\bx)= \rho(\bS^{-t}(\bx))$ with volume preservation $\mathrm{div} \,\bF = 0$.
Eigenfunctions of the Perron-Frobenius operator satisfy 
\begin{equation}
P^t \nu(\bx) = J^{-1}(\bx) \nu(\bS^{-1}(\bx)) = e^{\lambda t}\nu(\bx).
\end{equation}
Of particular interest is the
{\it physical} invariant measure $\mu^{*}(B) = \mu^{*}(\bS^{-t}(B))$ for all sets $B\in\mathfrak{B}$, which is stationary under the evolution of the flow.
The associated  invariant density 
$P\rho^{*}(\bx) = \rho^{*}(\bx)$
corresponds to an eigenfunction at eigenvalue $1$, which describes the asymptotic behavior of the underlying dynamics. Control is often designed to alter the observed invariant measure or density.
Spectral properties of the Perron-Frobenius operator are related, e.g., to almost-invariant sets, meta-stable states, mixing properties, decay of correlations~\cite{Gaspard1995pre,dellnitz1999jna,dellnitz2000nonl,Froyland2009pd}. 

The infinitesimal generator of the semigroup of Perron-Frobenius operators $\{ P^t \}_{t\geq 0}$ is given by the Liouville operator $L_P$~\cite{liouville1838jmpa,gaspard2005chaos}:
\begin{equation}\label{Eqn:LiouvilleEquation}
\frac{\partial }{\partial t} {\rho}(t,\bx) = - \nabla \cdot\left({\bf F}(\bx)\,{\rho}(t,\bx)\right) = L_P[{\rho}](\bx),\quad {\rho}_0:={\rho}(0,\bx),
\end{equation}
for continuously differentiable $\rho$ with compact support and appropriate boundary conditions.
The Liouville equation~\eqref{Eqn:LiouvilleEquation} describes how the flow transports densities in phase space and has a very intuitive interpretation as the conservation of probability or mass in terms of trajectories in the phase space. Alternatively, the evolution of the density may interpreted as propagated uncertainty of an initial state.
This is a first-order PDE that can be solved with the method of characteristics~\cite{zwillinger1989handbook,dutta2011jgcd}. 
The invariant density satisfies  $L_P[{\rho}^{*}](\bx)=0$ and is an eigenfunction at eigenvalue $0$.

\begin{figure}[tb]
	\centering
	\includegraphics[width=\linewidth]{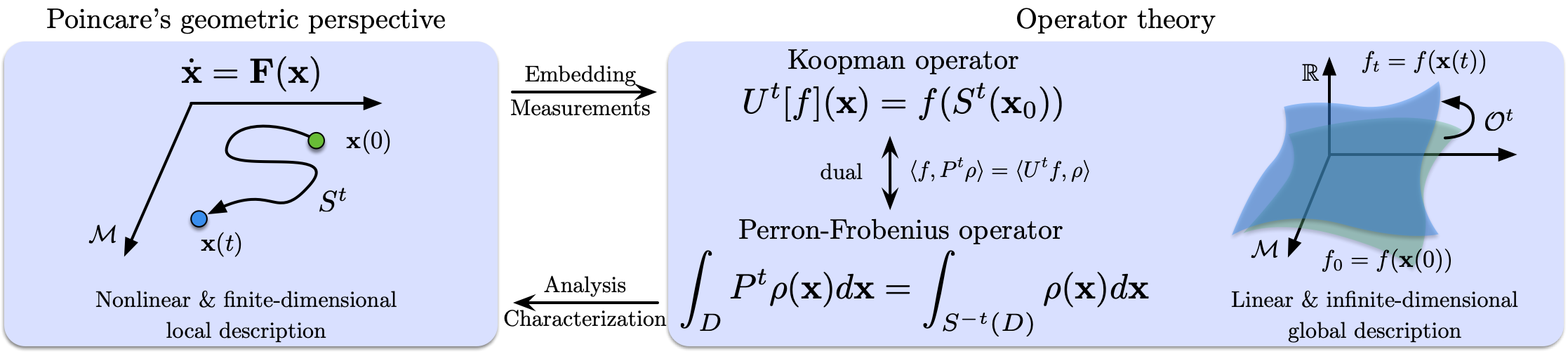}
	\caption{Poincar\'e's local phase space characterization versus the global operator-theoretic perspective.}
	\label{Fig:Geometry_vs_Operator}
\end{figure}

The aim of this chapter is to provide an overview of major developments and advances in data-driven control using operator-theoretic methods. 
The chapter is organized as follows:
In Sec.~\ref{Sec:ProblemFormulation}, the general control problem for a nonlinear system is formulated from an optimal control perspective. 
Control-dependent Koopman and Perron-Frobenius operators are introduced in Sec.~\ref{Sec:TransferOperators}. 
The main objective of the operator-theoretic approach is to find a linear representation of the underlying dynamical system. Data-driven system identification methods for the controlled nonlinear system are summarized in Sec.~\ref{Sec:SystemIdentification}. 
Important aspects in control theory such as observability, controllability, state estimation, and control design in the operator-theoretic framework are discussed in Sec.~\ref{Sec:Control}.
The chapter is concluded in Sec.~\ref{Sec:Conclusions} with a discussion on major challenges, open problems, and possible future directions of transfer operators approximations and their application for control.

%
%
\section{Control problem formulation}
\label{Sec:ProblemFormulation}
Reformulating strongly nonlinear dynamics in a linear framework via transfer operators is appealing as it enables the application of powerful optimal and robust estimation and control techniques available for linear systems~\cite{sp:book,dp:book,stengel2012book}.  
The formulation of an optimal control problem~\cite{stengel2012book} appears in many applications, such as trajectory control in robotics or boundary layer stabilization in fluids, and has been considered widely in the context of Koopman operators for control~\cite{Brunton2016plosone,Korda2016arxiv,Kaiser2017arxiv,KaKuBr2018arxiv}.

Optimization-based approaches for control provide a highly suitable framework for the control of nonlinear system, e.g. including constraints and flexible cost formulations. The optimal control problem can be stated as follows:
\begin{subequations}\label{Eqn:CostFunction}
\begin{align}
\min\limits_{\bu\in\mathbb{U}} & \int\limits_{t_0}^{t_f}\, l[\by(t),\bu(t),t] \,  d t +l_f[\by(t_f),t_f] 
	\end{align}
\end{subequations}
subject to
\begin{subequations}\label{Eqn:NonlinearSystemWithControl}
	\begin{align}
	\frac{d}{dt} \bx(t) &= {\mathbf{F}}(\bx,\bu),\quad \bx(0) = \bx_0\\
	\by(t) &= \bH (\bx,\bu)
	\end{align}
\end{subequations}
and possibly additional constraints on states, measurements, control inputs or time. 
We consider the nonlinear system~\eqref{Eqn:Dynamics} extended to include an external input ${\bu\in\mathbb{U}\subset\mathbb{R}^{q}}$, where $\mathbb{U}$ is the space of admissible control functions (or sequences in discrete time), and assume continuously differentiable state dynamics 
${\bf F}:\mathbb{X}\times\mathbb{U}\rightarrow \mathbb{X}$.
The state $\bx$ may not be fully accessible and instead a limited set of output measurements $\by\in\mathbb{Y}\subset\mathbb{R}^p$ may be collected, which are prescribed by the measurement function $\bH:\mathbb{X}\times \mathbb{U}\rightarrow\mathbb{Y}$.
The nature of the solution of the optimization problem is determined by the choice of the terminal state cost $l_f[\cdot]$ and the running cost $l[\cdot]$.
The objective is to determine a control law or policy that minimizes the cost functional~\eqref{Eqn:CostFunction}. 

A nonlinear optimal control formulation for the system in~\eqref{Eqn:NonlinearSystemWithControl} can be established using dynamic programming~\cite{Bertsekas2005book}. 
This relies on Bellman's principle of optimality~\cite{bellman1957book}
and leads to the Hamilton-Jacobi-Bellman (HJB) equation~\cite{bellman1964book}, a nonlinear PDE for the globally optimal solution.
Solving this nonlinear PDE is computationally challenging and as a result, a variational argument and Pontryagin's maximum principle~\cite{Pontryagin2062interscience} is often used instead, leading to a set of coupled ordinary differential equations (ODEs), the Euler-Lagrange equations.
While this two-point boundary value problem is solvable for higher-dimensional problems in contrast to the HJB equation, it is still too computationally demanding in real-time applications for many nonlinear systems, e.g. in autonomous flight~\cite{tang2018ar}.
Moreover, for high-dimensional systems, such as fluid flows, expensive direct and adjoint simulations render this approach infeasible, instead motivating the use of reduced-order models. 

Controlling nonlinear, and possibly high-dimensional, systems is generally computationally demanding and performance and robustness guarantees exist only for certain classes of dynamical systems. 
However, the control problem above simplifies considerably for quadratic cost functions and linear dynamics of the form
\begin{equation}\label{Eqn:LinearSystem}
\frac{\mathrm d}{\mathrm dt}\bx = \bA \bx + \bB\bu,\quad \bx(0) = \bx_0, 
\end{equation} 
where the matrix $\bA$ may be obtained by a suitable linearization of the nonlinear dynamics~\eqref{Eqn:NonlinearSystemWithControl} around an equilibrium point or operating condition.
Considering linear dynamics, the optimization problem can be simplified tremendously and becomes solvable even for high-dimensional systems. 
In the simplest case, without any constraints, this reduces to solving an algebraic Riccati equation (ARE) and yields the linear quadratic regulator (LQR)~\cite{stengel2012book}. 

Many problems may not permit a fully linear representation or the approximation may only be valid close to the linearization point. 
Moreover, models estimated from data may quickly become invalid, either due to the application of control or changing system parameters and conditions. 
For these cases, there exist adaptive and nonlinear variants of the formulation above. 
For instance, control-affine, nonlinear systems, whose governing equations can be factored into a linear-like structure, permit a state-dependent transition matrix $\bA(\bx)$ and actuation matrix $\bB(\bx)$ so that the ARE may be solved point-wise as a state-dependent Riccati equation (SDRE)~\cite{pearson1962ije}. The SDRE generalizes LQR for nonlinear systems, retaining a simple implementation and often yielding near optimal controllers~\cite{clautier1996proc}. 
Alternatively, multiple model systems can be considered, which usually consist of a large set of locally valid models, for which control is then determined either using a specific model selected based on some metric or by averaging/interpolating the control action from the model set~\cite{murray1997book}.

A widely adopted adaptive variant is model predictive control (MPC)~\cite{garcia1989model,morari1999model,allgower2004nonlinear}, which solves the optimal control problem over a receding horizon subject to the modeled dynamics and system constraints.
The receding horizon problem is formulated as an open-loop optimization over a finite time horizon. At each time step and given the current measurement, a sequence of future control inputs minimizing the cost $J$ over the time horizon is determined.
The first control value of this sequence is then applied, and the optimization is re-initialized and repeated at each subsequent time step. This
results in an implicit feedback control law $\bC(\bx_j) := \bu_{j+1}(\bx_j)$, where $\bu_{j+1}$  is the first input in the optimized actuation sequence starting at the initial condition $\bx_j:=\bx(t_j)$.
MPC is particularly ubiquitous in industrial applications including process industries~\cite{mayne2014automatica} and aerospace~\cite{eren2017jgcd}, as it enables more general formulations of control objectives and the control of strongly nonlinear systems with constraints, which are difficult to handle using traditional linear control approaches.

Linear models arising as approximations of the Koopman or Perron-Frobenius operators may be combined with any of these approaches. The potential benefit is the increased accuracy and validity as compared with linear models based on a local linearization and a reduced computational burden compared with schemes based on nonlinear models. Moreover, many data-driven approximation algorithms for these operators are easy to implement and efficient to compute as these are founded on standard linear algebra techniques.


%
%
\section{Control-oriented Transfer Operators}
\label{Sec:TransferOperators}
Controlled dynamical systems have been increasingly considered in the operator-theoretic framework, particularly for state estimation, to disambiguate dynamics from the effect of control, and for related problems such as optimized sensor and actuator placement.
While the Liouville equation has been examined for control purposes for a longer time~\cite{brockett2007ams,brockett2012chapter,roy2017jsc}, and interestingly also in the context of artificial intelligence~\cite{kwee2001ab}, the control-oriented formulations of the Koopman and Perron-Frobenius operators have gained traction only recently~\cite{Proctor2016arxiv,williams2016ifac,froyland2016siam,das2018arxiv}.


\subsection{Non-affine control systems}
We consider a non-affine control  dynamical system and  assume access to the full state $\by=\bx$:
\begin{equation}\label{Eqn:NonaffineControlSystem}
\frac{d}{dt} \bx(t) = {\mathbf{F}}(\bx,\bu),\quad \bx(0) = \bx_0.
\end{equation}
The flow associated with~\eqref{Eqn:NonaffineControlSystem} is referred to as the {\it control flow} and is denoted by $\bS^t(\bx,\bu)$.  Skew-product flows, such as $\bS^t(\bx,\bu)$, arise in topological dynamics to study non-autonomous systems, e.g.\ with explicit time dependency.
Thus, we consider the action of an operator in the extended state space $\mathbb{X}\times\mathbb{U}$.

The Koopman operator is defined as acting on the extended state:
\begin{equation}\label{NonaffineKoopmanControl}
U^t f(\bx,\bu) = f(S^t(\bx,\bu),\bu).
\end{equation}
Here the inputs $\bu$ evolve dynamically, e.g.\ by a prescribed exogenous behavior or a given state-dependent feedback control law $\bu = \bK(\bx)$.
If the inputs themselves are not evolving dynamically, this would reduce to $U^t f(\bx,\bu) = f(S^t(\bx,\bu),0)$, 
in which case the inputs parameterize the dynamics~\cite{Proctor2016arxiv}.
The spectral properties of $U$ should then contain information about the unforced system with $\bu={\bf 0}$.
Koopman eigenfunctions associated with~\eqref{NonaffineKoopmanControl} satisfy
\begin{equation}
U^t \phi(\bx,\bu) = e^{\lambda t} \phi(\bx,\bu).
\end{equation}
Assuming smooth dynamics and observables, the generator equation for the Koopman family is:
\begin{equation}
\frac{\partial}{\partial t} f(t,\bx,\bu) = \bF(\bx,\bu) \cdot\nabla_{\bx} f(t,\bx,\bu) + \dot{\bu} \cdot\nabla_{\bu} f(t,\bx,\bu). 
\end{equation}
As above, smooth Koopman eigenfunctions can be considered eigenfunctions of the infinitesimal generator satisfying
\begin{equation}
\frac{d}{dt} \phi(\bx,\bu) = \bF(\bx,\bu) \cdot\nabla_{\bx} \phi(\bx,\bu) + \dot{\bu} \cdot\nabla_{\bu} \phi(\bx,\bu)= \lambda \phi(\bx,\bu). 
\end{equation}

Representing the system in terms of a finite set of observables or eigenfunctions generally requires a reformulation of the cost functional~\eqref{Eqn:CostFunction} so that
\begin{equation}\label{Eqn:CostFunctionKO}
	J =  \int\limits_{t_0}^{t_f}\, l^K[\bz(t),\bu(t),t] \,  d t +l_f^K[\bz(t_f),t_f],
\end{equation}
where $\bz(t)= {\bf f}(\bx(t),\bu(t))$ describes the nonlinear transformation of the state $\bx$ through a vector-valued observable ${\bf f}(\bx,\bu)=[f_1(\bx,\bu),\ldots,f_d(\bx,\bu)]^T$.
It may be of interest to directly control specific observables, e.g. Koopman eigenfunctions, that are associated with a particular physical behavior. 
If the state $\bx$ is included as an observable, it is possible to transform the cost functional~\eqref{Eqn:CostFunctionKO} into~\eqref{Eqn:CostFunction} by modifying $l^K$ and $l_f^K$ accordingly. 
However, there does not exist a finite-dimensional Koopman-invariant subspace explicitly including the state, that is topologically conjugate to a system representing multiple fixed points, limit cycles, or more complicated structures~\cite{Brunton2016plosone}.
Nevertheless, it may be possible to obtain a linearization in the entire basin of attraction for a single fixed point or periodic orbit~\cite{Williams2015jnls,Lan2013physd}.
Further, Koopman eigenfunctions may be considered as observables and the state may be recovered via inversion, e.g. approximated from data using multidimensional scaling~\cite{kawahara2016nips}.


Analogously, we can consider the Perron-Frobenius operator acting on densities
\begin{equation}
	P^t \rho(\bx,\bu) = J^{-1}(\bx)\rho(S^{-t}(\bx,\bu),\bu),
\end{equation}
for which eigenfunctions satisfy
\begin{equation}
	P^t \nu(\bx,\bu) = e^{\lambda t} \nu(\bx,\bu). 
\end{equation}
The  control-oriented Liouville equation describing the transport of densities under the influence of an exogenous input is given by: 
\begin{equation}
	\frac{\partial }{\partial t} \rho(t,\bx,\bu) = - \nabla_{\bx}\cdot \left( \bF(\bx,\bu) \rho(t,\bx,\bu)\right) - \nabla_{\bu}\cdot \left( \dot{\bu} \rho(t,\bx,\bu)\right),\quad \rho_0(\bx,\bu) :=\rho(0,\bx,\bu),
\end{equation}
where the  initial condition is, e.g., $\rho_0(\bx,\bu) = \delta(\bx'-\bx)\delta(\bu'-\bu)$, a point-mass at $\bx'$ and $\bu'$ for deterministic dynamics.
This assumes the conservation of probability in the state-action space. The first term on the right-hand side describes the local change in density due to the flow, while the second term describes the local change due to exogenous inputs.

We are often interested in changing the long-term density $\rho(t,\bx)$.
The objective is then to determine inputs $\bu$ so that $\rho(t,\bx)$ becomes close  to a desired target density $\rho^{T}(\bx)$ over some time horizon $[0,t_f]$ or in the limit $t_f\rightarrow\infty$, which can be interpreted as steering a density of initial conditions or particles governed by the underlying dynamical system~\eqref{Eqn:NonaffineControlSystem}.
Assuming further that the inputs do not evolve dynamically, i.e.\ $\dot{\bf u} = 0$, we have:
\begin{equation}
\frac{\partial }{\partial t} \rho(t,\bx) = - \nabla_{\bx}\cdot \left( \bF(\bx,\bu) \rho(t,\bx)\right),
\end{equation}
which is considered in the seminal works of R. Brockett~\cite{brockett2007ams,brockett2012chapter}, particularly in the context of ensemble control. 
The finite-time control objective may be formulated as
\begin{equation}\label{Eqn:CostFunctionPFO}
J = \int_{t_0}^{t_f}\,\int_{\mathbb{X}}\left[\rho(t,\bx) - \rho^{T}(\bx)\right]l^P_{x}(\bx) d\bx + l^P_{u}[\bu(t),t]dt
\end{equation}
measuring the weighted deviation from the desired density and penalizing control input and time.
Ideally, the controlled system has a unique equilibrium density $\rho^{T}(\bx)$, corresponding to the sole eigenfunction of the Perron-Frobenius operator with eigenvalue $1$, with all others decaying exponentially.
In~\cite{brockett2007ams}, a general cost functional of the following form is proposed:
\begin{equation}\label{Eqn:CostFunctionPFO-1}
J = \int_{t_0}^{t_f}\int_{\mathbb{X}}\rho(t,\bx)l(\bx,\bu)d\bx dt+ \int_{\mathbb{X}} \left( \frac{\partial \bu}{\partial \bx}\right)^2d\bx + \int_{t_0}^{t_f}\left(\frac{\partial \bu}{\partial t}\right)^2 dt,
\end{equation}
where the first term evaluates the average performance of the system at the current time, the second term penalizes high gain feedback, and the third term may promote linearity of the control law.  

\begin{figure}[tb]
	\centering
	\includegraphics[width=\textwidth]{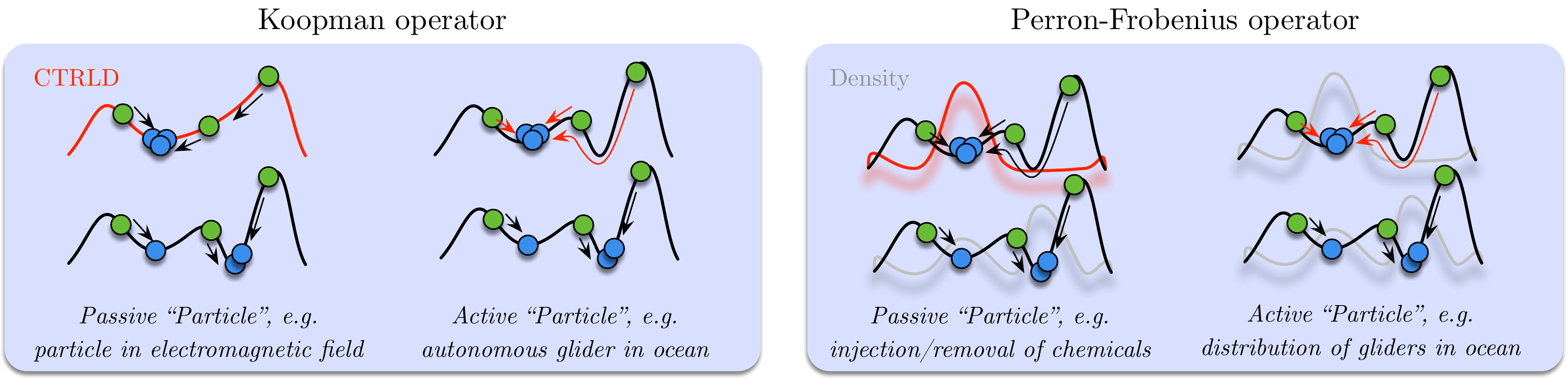}
	\caption{Controlling actively the observable or density function versus the trajectory in the operator-theoretic framework.}
	\label{Fig:ControlPerspectives}
\end{figure}

In general, different control perspectives can be assumed: (1) actively modifying the observable or density function, or (2) directly controlling the trajectory informed by the evolution dynamics of the observable or density function, e.g. ensemble control in a drift field.
An illustrative schematic is provided in Fig.~\ref{Fig:ControlPerspectives} using a particle in a field analogy.

%

\subsection{Affine control systems}
A control-affine dynamical system is given by
\begin{equation}\label{Eqn:NonlinearSystemWithControlAffine}
	\frac{d}{dt} \bx(t) = {\mathbf{F}}(\bx)+{\mathbf{G}}(\bx)\bu,\quad \bx(0) = \bx_0,
\end{equation}
where ${\mathbf{G}}(\bx)\bu:=\sum_i g_i({\bf x})u_i$ represents mixed terms in ${\bf x}$ and ${\bf u}$.  
Above we sought  a representation on the extended state $(\bx,\bu)$, in which case the spectral properties of the associated operators will contain knowledge on the unforced system with $\bu={\bf 0}$. Alternatively, we can consider operators associated with the unforced dynamics $\bF$ and how these are affected by the external control input, rendering the representation non-autonomous.
Then, additional terms will arise and inform how control affects observables, densities, and eigenfunctions:
\begin{equation}
\frac{\partial}{\partial t} f(t,\bx) = (\bF(\bx)+\bG(\bx)\bu)\cdot\nabla f(t,\bx)   =\bF(\bx)\cdot\nabla f(t,\bx) + \bG(\bx)\bu\cdot\nabla f(t,\bx),
\end{equation}
where $\nabla:=\nabla_{\bx}$.
We assume Koopman eigenfunctions are associated with the unforced dynamics $\bF(\bx)$.
Then, the smooth Koopman eigenfunctions satisfy
\begin{equation}
	\frac{d}{dt} \phi(\bx) = \lambda \phi(\bx) + \nabla \phi(\bx)\cdot \bG(\bx)\bu,
\end{equation}
which renders the dynamics governing the eigenfunctions generally nonlinear.
Only in the special case, where $\bG(\bx)$ is a linear function of the state $\bx$, this remains linear. The equation becomes bilinear for quadratic $\bG(\bx)$, for which a rich control literature exists~\cite{elliott2009bilinear}.

The Liouville equation for the density becomes:
\begin{equation}\label{Eqn:AffineSystem:ControlledLiouville}
\frac{\partial }{\partial t} \rho(t,\bx) = - \nabla\cdot \left( (\bF(\bx)+\bG(\bx)\bu) \rho(t,\bx)\right) 
= - \nabla\cdot \left( (\bF(\bx) \rho(t,\bx)\right) - \nabla\cdot \left(\bG(\bx)\rho(t,\bx)\right) \bu
\end{equation}
and correspondingly eigenfunctions satisfy:
\begin{equation}
\frac{d}{dt} \nu(\bx) = \lambda \nu(\bx) - \nabla\cdot(\bG(\bx)\nu(\bx))\bu. 
\end{equation}
The controlled Liouville equation has also been examined in other contexts, where a scalar observable is advected with the flow, e.g.\ concentration of chemical components or temperature distribution, in cases where diffusion is assumed negligible.
It may be possible that the controller directly injects or removes content to modify the density, e.g. 
for carbon removal strategies in the atmosphere or oil spill mitigation in the ocean.
Then, \eqref{Eqn:AffineSystem:ControlledLiouville} is modified to 
 \begin{equation}
 \frac{\partial }{\partial t} \rho(t,\bx) = - \nabla\cdot \left( (\bF(\bx) \rho(t,\bx)\right) + \sum_{i=1}^{p}\, g_i(\bx) u_i,
 \end{equation}
 where $g_i(\bx)$ is a set of functions prescribing the spatial impact of the actuators and the vector field $\bF$ remains constant in time.





%
%
\section{System identification}
\label{Sec:SystemIdentification}
In this section, we outline data-driven strategies to learn control-oriented representations of dynamical systems in the operator-theoretic framework.
Most approaches rely on dynamic mode decomposition for the Koopman operator and the classical Ulam's method for the Perron-Frobenius operator.
Emerging techniques such as machine learning and feature engineering aim to provide generalizations for systems of higher complexity, while exploiting and promoting sparsity often targets interpretability and efficient computation.
A general overview of algorithms for data-driven approximation of the Koopman and the Perron-Frobenius operators have been provided in~\cite{klus2017data}, which will not be repeated here. 
Here, we discuss approaches that specifically include the effect of control or external parameters. 


\subsection{Koopman operator}
\label{Sec:SystemIdentification_KO}


System identification for  the Koopman operator can be broadly divided into three main directions: (1) linear models as approximations to the Koopman operator describing the evolution of a, usually large, set of observable functions, (2) (bi)linear models for Koopman eigenfunctions, which constitute a Koopman-invariant subspace, and (3) linear models parametrized by the control input. 
We assume a collection of data tuples $(\bx_k,\bx_k',\bu_k)$, where the state $\bx_k:=\bx(t_k)\in\mathbb{R}^n$, the time-shifted state $\bx_k':=\bx(t_k+\Delta t)\in\mathbb{R}^n$, and inputs $\bu_k:=\bu(t_k)\in\mathbb{R}^q$ are sampled at discrete times $t_k$, $k=1,\ldots,m$. 
We define a vector-valued observable function ${\bf f}(\bx) = [f_1(\bx),f_2(\bx),\ldots,f_d(\bx)]^T$ with ${\bf f}:\mathbb{R}^n\rightarrow\mathbb{R}^d$ evaluated on data, which is generally a nonlinear transformation of the state so that $\bz = {\bf f}(\bx)$ embeds the state in a feature space of possibly higher dimension, i.e. $d\gg n$.
\newline

\textbf{Dynamic mode decomposition with control:}  
Building on DMD, the control-oriented extension {\it DMD with control} (DMDc) was proposed to disambiguate the effect of control from the natural dynamics in high-dimensional systems~\cite{Proctor2016siads}. 
Assuming the availability of full-state measurements, observables are taken as linear measurements of the state, i.e. ${\bf f} = \bx$.
Then, the dynamical behavior of coordinates $\bz$ are modeled as $\dot{\bz}(t) = \bA\bz(t) + \bB\bu(t)$ or 
\begin{equation}
{\bz}_{k+1} = \bA\bz_{k} + \bB\bu_{k}
\end{equation} 
in discrete time, which makes standard linear control techniques readily applicable.
The dimensionality of the problem is determined by the dimension of the state. However, many systems evolve on a low-dimensional attractor so that the singular value decomposition (SVD) can be used to solve the problem efficiently. Further, DMDc has been combined with compressed DMD~\cite{Brunton2015jcd} for compressive system identification~\cite{bai2017aiaa}.
While DMDc is feasible for high-dimensional systems with $n\ll m$,  such as fluids that may have millions or billions of degrees of freedom, linear measurements are known to not be sufficiently rich to capture strongly nonlinear behavior, motivating the use of nonlinear observables.
Nevertheless, recent work~\cite{Kaiser2018prsa} indicates that even models with relatively little predictive power may still be sufficient for model predictive control schemes (see also Sec.~\ref{Sec:ControlExamples}). 
\newline

\textbf{Extended DMD with control:} Extended DMD (eDMD) has been proposed as a generalization of DMD employing a dictionary of nonlinear functions as observables~\cite{Williams2015jnls}. DMD is then equivalent to eDMD with the special set of observables corresponding to linear monomials, ${\bf f}(\bx) = [x_1,x_2,\ldots,x_n]^T$. 
In domains where such a local linearization would be insufficient, more complex basis functions, such as high-order monomials, Hermite polynomials, or radial basis functions, can be employed. 
The number of observables typically exceeds the state dimension, i.e.\ $d\gg n$; however, the algorithm can be reformulated to scale with the dimension of samples instead of the dimension of features/states.
eDMD has also been extended for nonlinear system identification, which we refer to as {\it eDMD with control} (eDMDc)~\cite{williams2016ifac}, and been combined with MPC~\cite{Korda2016arxiv,Korda2018arxiv}. 
By transforming the system so that the dynamics evolve linearly, linear MPC methods become readily applicable, which require solving a convex quadratic program instead of a non-convex program as in nonlinear MPC.  In~\cite{Korda2016arxiv} the MPC optimization is further reformulated to scale with the dimension of the state $n$ and not with the number of observables $d$, which is advantageous if $d \gg n$. 
However, the quality of the resulting model and its prediction is highly sensitive to the choice of function basis. 
\newline

\textbf{Time-delay coordinates and DMDc:} Often we do not have access to the full state, but instead can only measure a single or few variables. Then time-delay coordinates, i.e. ${\bf f}(\bx(t)) = [x(t),x(t-1\tau),x(t-2\tau),\ldots,x(t-(d-1)\tau)]$ with time delay $\tau$, have been proven to be  useful for analyzing and modeling dynamics.
The Taken's embedding theorem~\cite{Takens1981lnm}
provides conditions under which the original system can be faithfully reconstructed in time-delay coordinates, providing a diffeomorphism between the two systems. 
Recently, time-delay coordinates have been shown to provide a Koopman-invariant subspace and used to model unforced systems with bimodal or bursting behavior as linear systems closed by an intrinsic forcing term~\cite{Brunton2017natcomm}. 
This formulation is achieved by performing DMD on the delay coordinates arranged in a Hankel matrix~\cite{Brunton2017natcomm,Arbabi2016arxiv}.  
Hankel matrix based system identification has been used for decades in the eigensystem realization algorithm (ERA)~\cite{ERA:1985}. 
Assuming sample pairs $(x_k,u_k)$ of a single long trajectory $x(t)$ generated by a continuous system with single control input $u(t)$, 
the vector-valued observable is given by $\bz_k = {\bf f}(\bx_k) = [x_k, x_{k-1}, x_{k-2},\ldots, x_{k-d_1+1}]$ where $f_{i}(x_{k})=x_{k-(i-1)}$ and $\bv_k =  [u_k, u_{k-1}, u_{k-2},\ldots, u_{k-d_2+1}]$.
The embedding dimensions $d_1=d_2=d$ are assumed equal for simplicity. 
The system is then given by
\begin{subequations}\label{Eqn:KM:TDC_MISO}
	\begin{align}
	{\bz}_{k+1} = \bA{\bz}_k + \bB \bv_k\\
	\bx_{k} = \begin{bmatrix}
	1 & 0 & \ldots & 0
	\end{bmatrix}{\bz}_k
	\end{align}
\end{subequations}
where the state is recovered from the first component of $\bz$.
The control matrix $\bB$ must satisfy a lower triangular structure to not violate causality.
If the number of time-delay coordinates is large, an SVD can be applied and the model is built on the first $r$ eigen-time delay coordinates. 
It can be advantageous to embed the full state $\bx$ if available instead of just a single component, as more information is used to improve prediction accuracy.
In Eq.~\eqref{Eqn:KM:TDC_MISO} the control history appears as a separate vector of control inputs; it can also be useful to augment the state with the past control values so that the current actuation value appears as single control input:
\begin{subequations}\label{Eqn:KM:TDC_SISO}
	\begin{align}
	\hat{\bz}_{k+1} = 
	\begin{bmatrix}
	{\bz}_k\\
	u_{k-1}\\
	\vdots\\
	u_{k-d+1}
	\end{bmatrix}_{k+1} &= 
	\left[\begin{array}{c|c}
	\bA & \bB_{[2,d]} \\
	\hline
	{\bf 0} & \begin{bmatrix}
	0 & {\bf 0}\\
	{\bf 0} & {\bf I}
	\end{bmatrix}
	\end{array}\right]
	\begin{bmatrix}
	{\bz}_k\\
	u_{k-1}\\
	\vdots\\
	u_{k-d+1}
	\end{bmatrix}_{k} + 
	\begin{bmatrix}
	\bB_{[1]} \\ 1 \\ 0 \\ \vdots \\0
	\end{bmatrix}
	u_k\\
	\bx_{k} &= \begin{bmatrix}
	1 & 0 & \ldots & 0
	\end{bmatrix}
	\hat{\bz}_k,
	\end{align}
\end{subequations}
where ${\bf I}$ denotes the $d-2\times d-2$ identity matrix and $\bB_{[a,b]}$ contain columns $a$ through $b$ of $\bB$. 
\newline

\textbf{Control in eigenfunction coordinates:} Eigenfunctions of the Koopman operator provide a natural set of intrinsic observables as these themselves behave linearly in time and correspond to global properties of the system. 
Particularly, control affine systems Eq.~\eqref{Eqn:NonlinearSystemWithControlAffine} have been studied in this context. 
The general approach relies on estimating eigenfunctions associated with the autonomous system and then identifying how these are affected by the control input.  
The observables are then given by a set of eigenfunctions $\bz = {\bf f}(\bx) = [\phi_1(\bx), \ldots, \phi_r(\bx)]$.  Their continuous-time dynamics are given by
\begin{equation}\label{Eqn:KM:KRONIC}
\dot{\bz} = (\bF(\bx)+\bG(\bx)\bu)\cdot\nabla_{\bx}{\bf f}(\bx) = \bA \bz + \bB(\bz) \bu,
\end{equation}
where $\bA := \bF(\bx)\cdot\nabla_{\bx}{\bf f}(\bx) = \bLambda := \mathrm{diag}(\lambda_1,\ldots,\lambda_r)$
and $\bB(\bz)  := \bG(\bx)\bu\cdot\nabla_{\bx}{\bf f}(\bx)$.
For $\bu = {\bf 0}$, this constitutes a Koopman-invariant subspace for the unforced system; otherwise, the control term $\bB(\bz)$ prescribes how these eigenfunctions are modified by $\bu$.
The advantage of this formulation is that the dimension of the system scales with the number of eigenfunctions, and we are often interested in a few dominant eigenfunctions associated with persistent dynamics. 
This formulation has been examined for state estimation and observer design~\cite{Surana2016cdc,surana2017cdc} and data-driven control~\cite{Kaiser2017arxiv}.
In general, the eigenfunctions can be identified using eDMD, kernel-DMD or other variants, and the model~\eqref{Eqn:KM:KRONIC} may be well represented as long as their span contains $\bF(\bx)$, $\bG(\bx)$, and $\bx$~\cite{surana2017cdc}.  
However, the eigenfunctions obtained from DMD/eDMD may be spurious, i.e.\ they do not behave linearly. 
For instance, noise can produce corrupt, physically irrelevant but energetically important, unstable eigenfunctions. 
An alternative~\cite{Kaiser2017arxiv} seeks a functional representation of the eigenfunctions in a dictionary of basis functions $\boldsymbol{ \Theta}(\bx) = [\theta_1(\bx),\theta_2(\bx),\ldots,\theta_p(\bx)]$ so that $\phi(\bx)\approx \sum_{i=1}^p\,\theta_i(\bx)\xi_i = \boldsymbol{\Theta}(\bx)\mathbf{\xi}$.
A sparsity constraint on the vector of coefficients $\xi\in\mathbb{R}^d$  can then be used to identify an analytical expression of the eigenfunction. This approach is restricted to smooth eigenfunctions in the point spectrum of the Koopman operator and results in an optimization problem to find a sparse vector in the nullspace of the matrix $\left(\dot{\bx}\cdot\nabla_{\bx} \boldsymbol{\Theta}(\bx) - \lambda \boldsymbol{\Theta}(\bx)\right)\mathbf{\xi} = 0$ for a given eigenvalue $\lambda$. 
This requires an accurate estimation of the time derivative of $\bx$ and is sensitive to noise. The control term $\bB(\bz)$ can either be determined from the gradient of the identified eigenfunctions for given $\bG(\bx)$ or identified directly~\cite{KaKuBr2018arxiv} for unknown $\bG(\bx)$ by similarly expanding it in terms of a dictionary, $\bG(\bx) \approx \sum_{i=1}^q \psi_i(\bx) \eta_i$.
The approach can be further extended to nonaffine control systems, for which control-dependent eigenfunctions may be considered, using a library on the extended state, i.e.\ $\boldsymbol{ \Theta}(\bx,\bu)$.
This data-driven characterization has been combined with SDRE~\cite{Kaiser2017arxiv} and MPC~\cite{KaKuBr2018arxiv} for controlling the underlying dynamical system. 
The cost functional generally needs to be formulated in eigenfunction coordinates as discussed in the previous section.
\newline

\textbf{Parameterized linear models: } Another class of linear parameter-varying models have been considered to disambiguate the approximation of the Koopman operator of the unforced system from the effect of exogenous forcing~\cite{williams2016ifac} and for model predictive control of switched systems~\cite{Peitz2017arxiv}. 
In particular, a set of linear models parametrized by the control input is considered as an approximation to the non-autonomous Koopman operator:
\begin{equation}\label{Eqn:KM:ParametrizedKoopman}
{\bz}_{k+1} = \bA (\bu) \bz_k.
\end{equation} 
For discrete-valued inputs $\bu_i\in\{\bu_1,\ldots,\bu_q \}$, a finite set of matrices defined as $\{\bA_i:=\bA(\bu_i)\}_{i=1}^q$ may be estimated, in the simplest case with DMD or a variant by determining these for each discrete control input separately.
This replaces the non-autonomous Koopman operator with a set of autonomous Koopman operators,  each element associated with a constant control input, as considered in~\cite{Proctor2017siads}.
A more robust way that generalizes to continuous input data has been proposed in~\cite{williams2016ifac}, where $\bA(\bu)$ is expanded in terms of basis functions $\psi_i:\mathbb{U}\rightarrow\mathbb{R}$ as in eDMD for the state variable, so that $\bA (\bu) = \sum_{k=1}^{d} \psi_k(\bu) \bA_k$. Then, for a given value of $\bu$, matrix $\bA$ is constructed using coefficients $\bA_k$ and its spectral properties can be examined. Since eDMD may be prone to overfitting, the estimation problem can be regularized using the group sparsity penalty~\cite{simon2012standardization}.   

Determining an optimal control sequence for Eq.~\eqref{Eqn:KM:ParametrizedKoopman} is a combinatorial problem, which can be efficiently solved for low-dimensional problems using iterative algorithms from dynamic programming~\cite{Bertsekas2005book}.
In~\cite{Peitz2017arxiv}, the MPC problem for the sequence of control inputs is transformed into a time switching optimization problem assuming a fixed sequence of consecutive discrete control inputs. The model can further be modified as a bilinear continuous, piecewise-smooth variant with constant system matrices, which does not suffer from the curse of dimensionality and allows for continuous control inputs via linear interpolation~\cite{peitz2018feedback}. 


\subsection{Perron-Frobenius operator}
\label{Sec:SystemIdentification_PFO}
A classical method to approximate the Perron-Frobenius operator is the Ulam-Galerkin method, a particular Galerkin method where the test functions are characteristic functions. Given an equipartition of the phase space into disjoint sets $\{B_1,\ldots,B_d\}$ and data, the action of the operator can be approximated as 
\begin{equation}\label{Eqn:TransitionMatrix}
\bP_{ij}^{\tau} = \frac{\mathrm{card}( \bx_k \vert \bx_k \in B_j \wedge \bF^{\tau}(\bx_k) \in B_i)}{\mathrm{card}(\bx_k \in B_j)},
\end{equation} 
where $\bF^{\tau}$ is the flow map and $\tau$ is the time duration.
The $\tau-$step stochastic transition matrix acts on a discrete probability vector $\bp=[p_1,\ldots,p_d]^T$ and satisfies the properties $\sum_{j=1^d}P_{ij}=1$ and  $0\leq P_{ij}\leq 1$. A widely used software package based on set-oriented methods is GAIO~\cite{dellnitz2001gaio}. 
Each element $B_i$ can be associated with a distinct symbol. The dynamics are then propagated on the partition and can be analyzed in terms of the sequence of symbols.
The approach has been generalized to PDEs~\cite{vaidya2009cdc,Kaiser2014jfm}, i.e. where the evolution of a pdf $p(\bu,t)$ of a spatiotemporal vector field $\bu(\bx,t)$ is of interest. After expanding the vector field $\bu(\bx,t) = \sum_i a_i(t)\bu_i(\bx)$, e.g. using POD, Ulam's method can be applied to the $a_i$'s.
Other data-driven methods to approximate the Perron-Frobenius operator include blind source separation, the variational approach for conformation dynamics (VAC)~\cite{noe2013variational} 
and DMD-based variants exploiting the duality between Koopman and Perron-Frobenius operators, e.g. naturally structured DMD (NSDMD)~\cite{huang2017arxiv} and constrained Ulam dynamic mode decomposition~\cite{goswami2018csl}. 

Data-driven control based on the Perron-Frobenius operator appears more challenging as it relies on good statistical estimates requiring large amounts of data.
In practice, methods often use a Monte-Carlo sampling to estimate the transition probabilities, which suffers from the curse of dimensionality. 
However, it is particularly appealing as it provides a natural framework to incorporate uncertainties and transients, allows the control of ensembles and mixing properties of the underlying dynamical system, and yields a more global controller.
\newline

\textbf{Parametrized models: }
Identifying a probabilistic model disambiguating the unforced dynamics and the effect of actuation is commonly realized by a parametrized representation.
Incorporating the effect of control can be achieved via a simple extension of Ulam's method by considering a discretized control input that parametrizes the transition matrix in~\eqref{Eqn:TransitionMatrix}: 
\begin{equation}\label{Eqn:ControlTransitionMatrix}
\bP_{ij}^{u_l} = \frac{\mathrm{card}( \bx_k \vert \bx_k \in B_j \wedge \bF^{\tau}(\bx_k,\bu_l) \in B_i)}{\mathrm{card}(\bx_k \in B_j)},
\end{equation} 
where each $\bP(\bu_l):=\bP_{ij}^{u_l}$ satisfies the normalization and positivity properties. 
The dynamics of the probability vector can then be expressed as
\begin{equation}\label{Eqn:DSDT-MarkovModel}
\bp_{k+1} = \bP(\bu_l) \bp_k.
\end{equation} 
Approximations for~\eqref{Eqn:ControlTransitionMatrix} have been obtained employing set-oriented methods~\cite{das2017acc} and from NSDMD and using Gaussian radial basis functions~\cite{das2018arxiv}, which are then used to derive optimal feedback control strategies.
The model~\eqref{Eqn:DSDT-MarkovModel} can be derived from a discretization of the Liouville equation in space, time, and control~\cite{kaiser2017tcfd} and represents a control-dependent Markov chain.
The corresponding control problem can be associated with a Markov decision process (MDP) and appears in discrete optimal control~\cite{Bertsekas2005book}. 
Therein, it is used as a model to predict the expectation of the value function in dynamic programming, without reference to the Liouville equation. 
Building on these ideas, cluster-dependent feedback control has been developed by examining the spectral properties of the control-dependent transition matrix family~\cite{kaiser2017tcfd}, where state discretization of a fluid from high-dimensional snapshot data is achieved using POD and k-means clustering.
\newline

\textbf{Multiplicative models: }
Alternatively, the effect of control may be encoded in a separate stochastic matrix $\bP^u$ approximating the action of a stochastic kernel: 
\begin{equation}
\bp_{k+1} = \bP \bP^u\bp_k,
\end{equation} 
which recovers the unforced dynamics when $\bP^u = {\bf I}$ is the identity matrix.
This model appears, e.g., in the context of optimizing mixing properties~\cite{froyland2016siam,froyland2016arxiv}, where $\bP$ represents advection via the Perron-Frobenius operator and  $\bP^u$ represents the discretized diffusion kernel to be optimized.
In~\cite{mehta2008tac}, limitations of nonlinear stabilization are examined when $\bP^u = \bP^{1+K}$ encodes a state feedback law. 
\newline



\textbf{Additive models: }It may also be possible to formulate additive models, which exhibit some similarity to traditional state-space models, of the form
\begin{equation}\label{Eqn:PFO_additive}
\bp_{k+1} = (\bP  + \bP^u)\bp_{k},
\end{equation}  
where additional constraints are necessary, so that $\bP  + \bP^u$ satisfies the positivity and conservation constraints from probability. The matrix $\bP^u$ may be interpreted as a disturbance to the unforced transition probabilities $\bP$. 
It may also be possible to represent $\bP^u$ as $\bP^c\bu_k$, where $\bP^c$ is a constant matrix and $\bu_k$ is the control input.
A similar formulation to~\eqref{Eqn:PFO_additive} incorporating control appears from a state-discretization of the controlled Liouville equation as outlined in~\cite{Kaiser2014jfm}, resulting in a discrete-state, continuous-time Markov model.
Variations of this type of model have been used to study optimal response of Markov chains~\cite{antown2018arxiv} and to control chaos~\cite{bollt2000ijbc}.
\newline

\textbf{Remark: }
Until recently, the focus on control in the  Perron-Frobenius and Liouville operator framework has been mostly analytic and theoretical. Further, most problems studied have been low-dimensional. 
On the other hand,  there exists an extensive literature on probabilistic models for control, 
which can benefit the operator-theoretic framework but which have not been fully explored yet.
Examples include Hidden Markov models, see~\cite{alpcan2006cdc,yu2007cdc} in the context of the Perron-Frobenius and Fokker-Planck operators, and partial-observation Markov decision processes (POMDP)~\cite{aastrom1965jmaa}.

\subsection{Numerical example}
\label{Sec:ControlExamples}
\begin{figure}[tb]
	\centering
	\begin{overpic}[width=\linewidth]{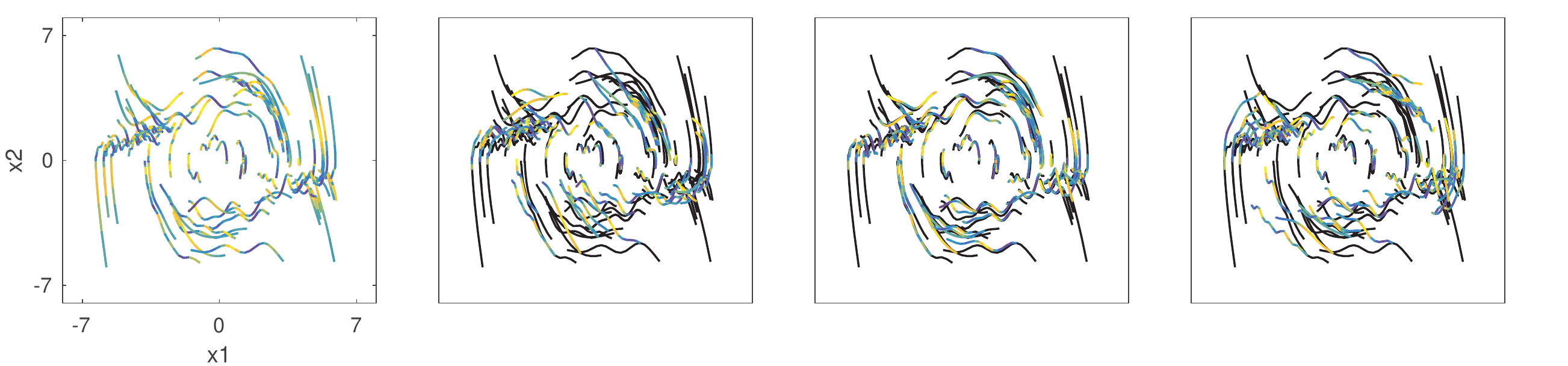}
		\put(-2,4){(a)}
		\put(5.5,5){\colorbox{white}{Validation data}}
		\put(33,24){\colorbox{white}{DMDc}}
		\put(57,24){\colorbox{white}{eDMDc}}
		\put(80,24){\colorbox{white}{delayDMDc}}
	\end{overpic}\\
	\begin{overpic}[width=\linewidth]{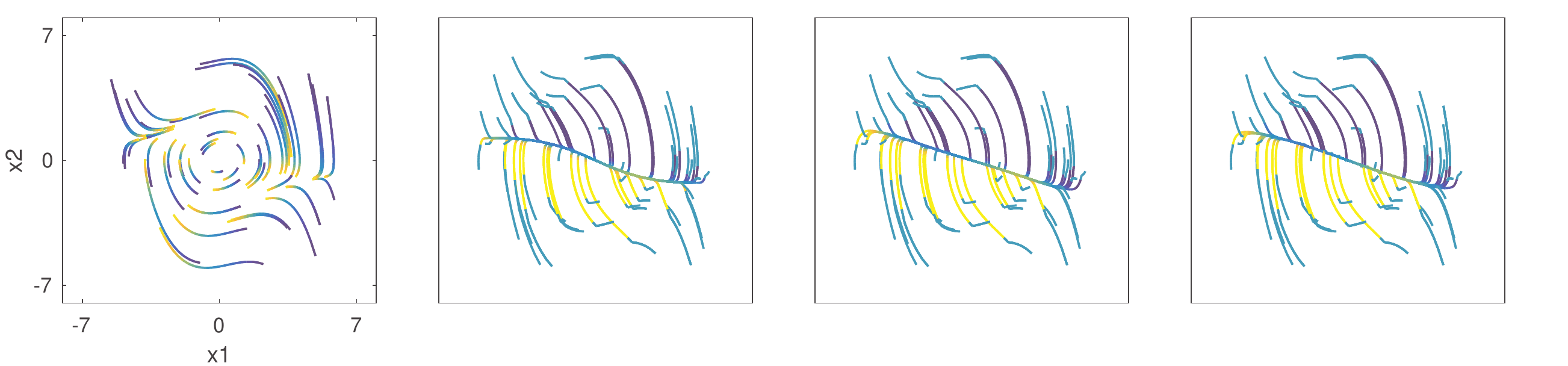}
		\put(-2,4){(b)}		
		\put(8,5){\colorbox{white}{Unforced}}
		\put(32,5){\colorbox{white}{Controlled}}
	\end{overpic}\\		
	\begin{overpic}[width=\linewidth]{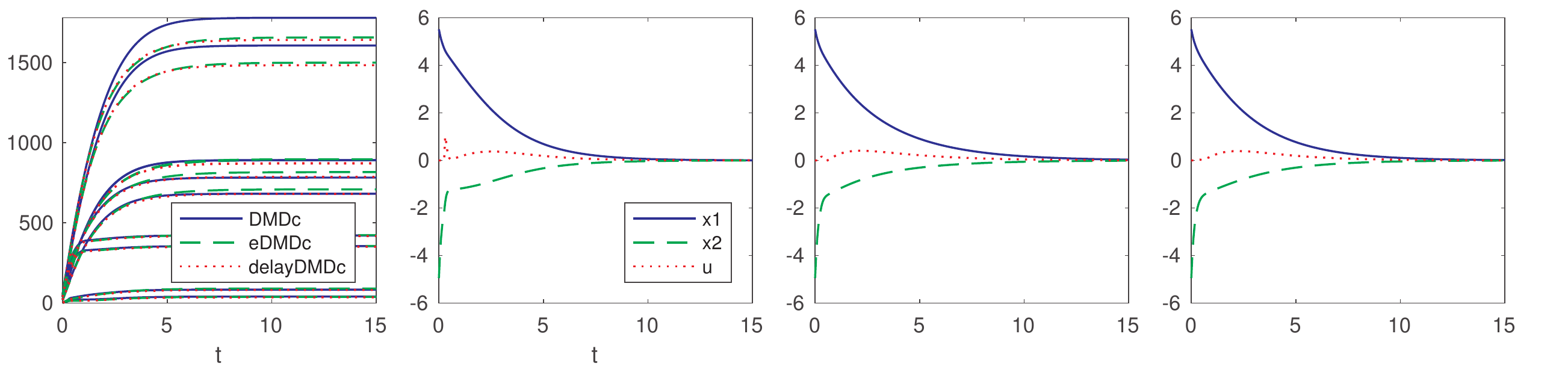}
		\put(-2,4){(c)}		
	\end{overpic}
	\\		
	\begin{overpic}[width=\linewidth]{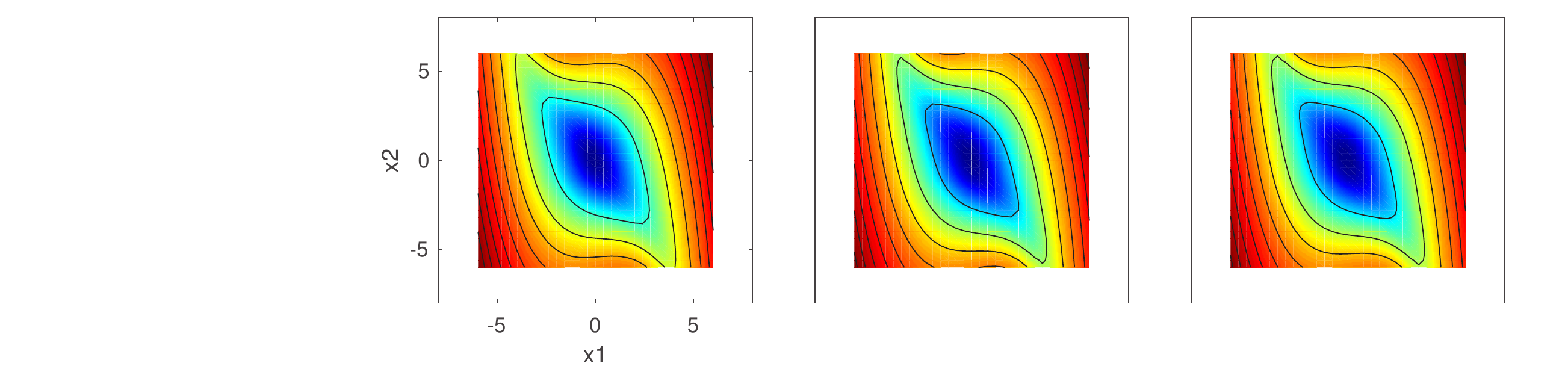}
		\put(-2,4){(d)}		
	\end{overpic}
	\caption{Van der Pol oscillator for $\mu = 0.2$: 
		(a) validation data and prediction (colored by control input),
		(b) unforced and forced phase plots with MPC for initial conditions of the validation dataset (colored by control input),
		(c) cumulative cost of subset of the trajectories and a specific example trajectory, and
		(d) color-coded initial conditions with respect to converged control cost. 
	}
	\label{Fig:MPC-Example}
\end{figure}
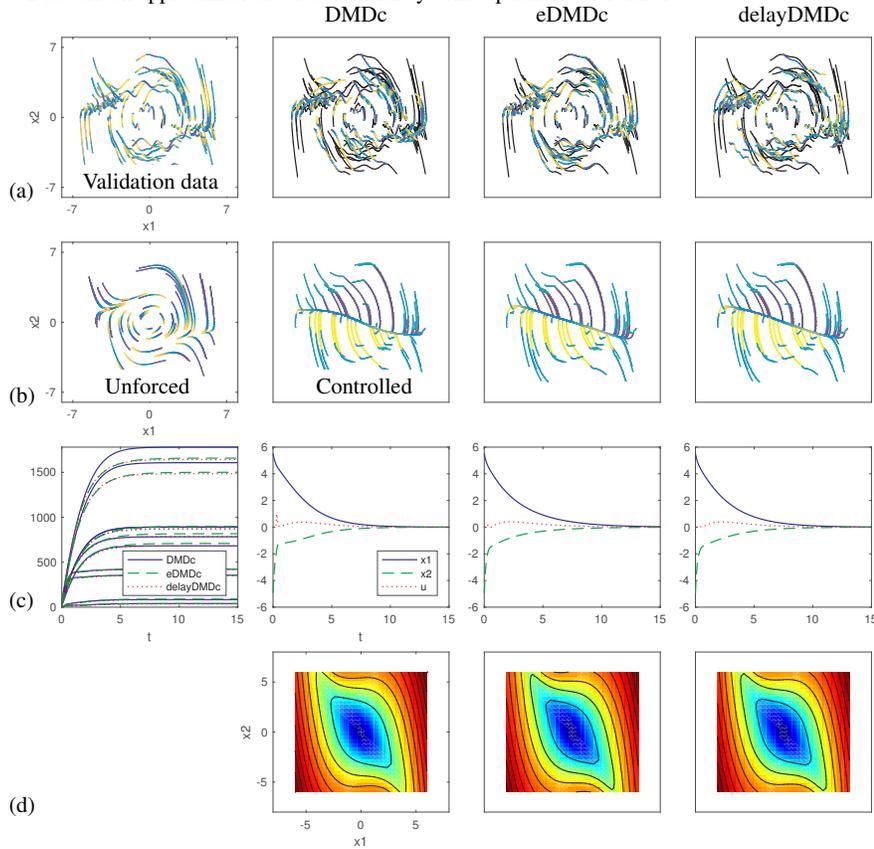

We compare the effectiveness of model predictive control for different DMDc-based models, which differ in their choice of dictionary functions. 
The van der Pol oscillator is considered as an illustrative example:
\begin{equation}
\frac{d^2 x}{d t^2} - \mu(1-x^2)\frac{d x}{d t}  + x = u
\end{equation}
with $\mu = 0.2$ and where $u$ is the external control input.
In the following, the state vector is defined as $\bx:=[x_1,x_2]^T=[x,dx/dt]^T$.
The training data consists of 200 trajectories in the box $[-6,\, 6]\times [-6,\, 6]$ integrated until $T=1$ with timestep $\Delta t=0.05$, i.e. 4000 sample points in total. The forcing $u(t) = 5 \sin(|\omega_1| t)\sin(|\omega_2| t)$, where $\omega_i\sim\mathcal{N}(0,10)$, is applied for the identification task.
The models under consideration include DMDc on the state $\bx$, eDMDc using monomials up to order $5$, and delayDMDc with embedding dimension $d=5$, as compared in Fig.~\ref{Fig:MPC-Example}. 
For the control optimization, weight matrices $Q = (\begin{smallmatrix}
1 & 0\\ 0 &1
\end{smallmatrix})$, $R_u = 0.1$, and $R_{\Delta u}=0.1$, and input constraints $-5<u<5$, $-50<\Delta u<50$ are employed. The prediction/control horizon is $T = 0.75= 15 \Delta t$. 

The models are evaluated on a test dataset that is different from the training set, visualized for a subset of trajectories in Fig.~\ref{Fig:MPC-Example}(a). 
The smallest error is achieved by eDMDc, followed by delayDMDc and DMDc, with the latter exhibiting a large deviation from the actual evolution.
MPC control results are displayed for $\mu = 0.2$ in Fig.~\ref{Fig:MPC-Example}(b)-(c).
Three main observations can be made: 
(1) Despite the prediction inaccuracies, all three models are able to successfully control the system;
(2) best performance is achieved with eDMDc and delayDMDc, while delayDMDc generally performs slightly better; 
(3) the farther away the initial condition is from the fixed point to be stabilized, the worse the performance for DMDc (compare results in Fig.~\ref{Fig:MPC-Example}(d)).
It is important to point out that the specific control performance is highly sensitive to the chosen prediction horizon, weights, and the forcing for the training data; however, delayDMDc appears to be most robust. 
Observed here and in previous work~\cite{Kaiser2018prsa}, superior prediction performance may be overrated in combination with predictive control schemes. The overall performance lies more in the robustness of MPC than in the specific model choice due to the repeated initialization with the updated measurement vector and the sufficiently short prediction horizon.

%
%
\section{Control theory in the operator-theoretic framework}
\label{Sec:Control}
In this section, we review control-theoretic advances using transfer operators and their connections to classical concepts in control such as stability analysis and observability.

\subsection{Stability Analysis}\label{Sec:StabilityAnalysis}

\begin{figure}[tb]
	\centering
	\includegraphics[width=\linewidth]{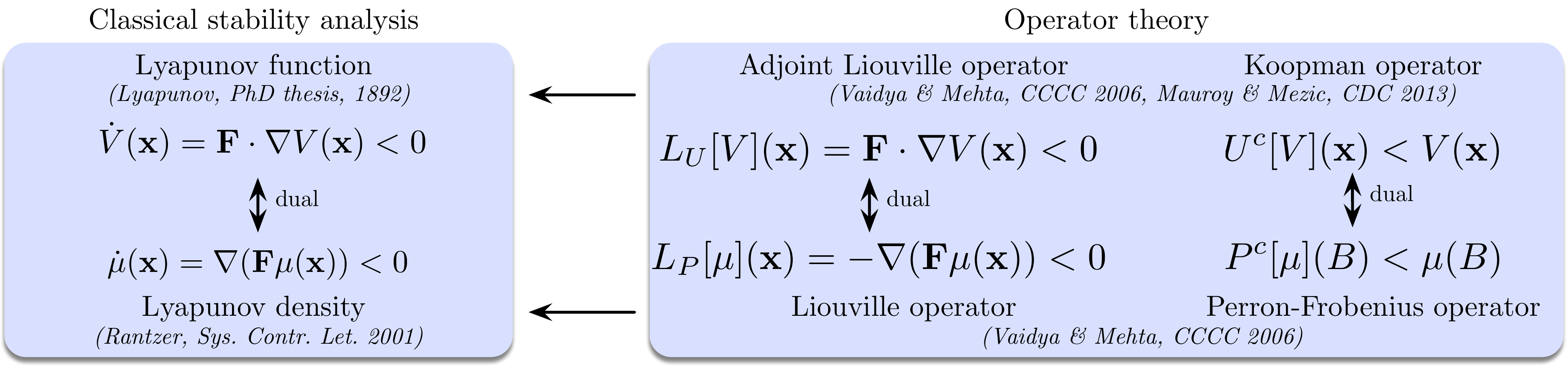}
	\caption{Classical stability analysis and within the operator-theoretic framework. The duality between the Lyapunov function and Lyapunov density~\cite{rantzer2001scl} can be derived from the duality of the adjoint transport operators. Operators $U^c$ and $P^c$ are restricted to the corresponding function space with support on $\mathbb{X}\backslash \mathbb{A}$, where $\mathbb{A}$ is the attractor.}
	\label{Fig:StabilityAnalysis}
\end{figure}

Analyzing the stability of dynamical systems is an important aspect of control theory, e.g. in order to develop stabilizing controllers. 
In particular, the Lyapunov function plays an essential role in the global stability analysis and control design of nonlinear systems. 
For stable, linear systems, the Lyapunov function is the positive solution of a corresponding Lyapunov equation. 
In contrast, finding and constructing a Lyapunov function poses a severe challenge for general nonlinear systems.
There has been important progress within the operator-theoretic context for nonlinear stability analysis as spectral properties of the operators can be associated with geometric properties such as Lyapunov functions and measures, contracting metrics, and isostables.

Some of this work relates back to a weaker notion of stability in terms of a density function as introduced by A. Rantzer~\cite{rantzer2001scl}.
In particular, he introduces the concept of asymptotic stability in an almost everywhere (a.e.) sense, which is guaranteed through the existence of a density function.
Around the same time, it was shown that Lyapunov functions can be constructed using a suitable choice of polynomials, resulting in a linear problem formulation~\cite{Parrilo2000phd}.
These ideas were later combined in a convex, linear formulation based on the discretization of the density function via polynomials, to jointly synthesize the density function and state feedback.
The crucial connection of this line of work to the Perron-Frobenius operator was made by U. Vaidya and P. Mehta~\cite{vaidya2006ccc} demonstrating that these ideas of duality and linearity can be expressed using methods from ergodic theory, and the Lyapunov measure represents the dual of the Lyapunov function, including its relationship to the Koopman operator. 
Spectral properties, specifically invariant sets, are then used to determine stability properties. The Lyapunov measure is closely related to the density function introduced in~\cite{rantzer2001scl} and captures the weaker a.e. notion of stability. 
Analogously to the Lyapunov equation, the solution to the Lyapunov measure equation provides necessary and sufficient conditions for a.e. stability of nonlinear systems. 
The Lyapunov equation and Lyapunov measure have been connected to the Koopman and Perron-Frobenius operator, respectively, e.g. providing an explicit formula for the Lyapunov measure in terms of the Perron-Frobenius operator, which generalizes stability analysis to invariant sets of nonlinear systems~\cite{vaidya2007acc}.
Set-oriented methods have been used for solving the Lyapunov measure equation leading to a linear program~\cite{vaidya2008tac}. 
A weaker notion of {\it coarse stability} has been further introduced to characterize the stability properties obtained from the finite-dimensional approximation.
A formulation for the infinitesimal generator of the Perron-Frobenius operator, the Liouville operator, has also been proposed and can be solved numerically using the finite-element method~\cite{rajaram2010jmaa}.
In particular, for a stable fixed point $\bx^{*}\in\mathbb{X}$, a Lyapunov measure (or density) $\mu$ satisfies a.e.
\begin{equation}
L_P [\mu] (\bx)<0\quad\forall \bx{\not =} \bx^{*}, 
\end{equation}
which corresponds to $\nabla (\bF \mu(\bx))>0$ as originally introduced in~\cite{rantzer2001scl}.
Thus, the Lyapunov measure (density) decays a.e. under the action of the Liouville operator and is related to its spectral properties.
These ideas have further been extended to stochastic systems~\cite{vaidya2015acc,vaidya2015cdc}, used to compute the domain of attraction~\cite{wang2010cdc}, and employed for control~\cite{vaidya2010ieee,raghunathan2014ieee}.

The Koopman operator has also been used to study nonlinear stability. 
Certain sets of points in phase space that exhibit the same asymptotic behavior are referred to as isostables and play an important role in characterizing transient behavior.
Isostables correspond to level sets of a particular Koopman eigenfunction and can be computed in the entire basin of attraction using, e.g., Laplace averages~\cite{mauroy2013pd}.
Building on~\cite{vaidya2006ccc}, it has been further shown~\cite{mauroy2013pd,mauroy2013cdc} that special Lyapunov functions from classical stability theory can be determined from Koopman eigenfunctions and a sufficient condition for global stability is the existence of stable eigenfunctions. 
As with the Lyapunov measure, it can be shown that for a stable fixed point there exists 
a special, nonnegative observable $V(\bx)$ that satisfies:
\begin{equation}
	L_U  [V](\bx)<0\quad\forall \bx{\not =} \bx^{*}. 
\end{equation}
Thus, the Lyapunov function decays everywhere under action of the adjoint Liouville operator and is related to its properties and those of the Koopman operator family.
In~\cite{mauroy2013cdc}, it has been shown that global stability of a fixed point can be established through the existence of a set of $C^1$ eigenfunctions of the Koopman operator associated with the eigenvalues of the Jacobian of the vector field, and that Koopman eigenfunctions can be used to define a Lyapunov function and contracting metrics. 
A numerical scheme based on a Taylor expansion is proposed to compute stability properties including the domain of attraction~\cite{mauroy2013cdc}.
These ideas have been further extended to examine global stability properties of hyperbolic fixed points and limit cycles and also for non-analytical eigenfunctions~\cite{Mauroy2016ieee}.

\subsection{Observability and controllability}
Characterizing a system's observability and controllability is crucial for designing sensor-based estimation and control strategies. 
Linear observability and controllability criteria have been extended to characterize nonlinear systems~\cite{Hermann1977ieeetac,khalil1996book}, e.g. local observability can be tested by constructing an observability matrix using Lie derivatives.
Balanced model reduction has also been extended to nonlinear systems~\cite{scherpen1993scl,lall2002ijrnc,zhang2002automatica}. 
Operator-theoretic approaches provide new opportunities for determining observability and controllability properties of high-dimensional nonlinear systems from data using linear techniques.

There has been considerable work on nonlinear observability analysis based on the Perron-Frobenius operator formulation for the Lyapunov measure equation. 
Analogously to the observability Gramian for linear systems, which can be obtained as the solution to a Lyapunov equation, the nonlinear equivalent can be obtained from the Lyapunov measure equation. 
In~\cite{vaidya2007cdc}, set-oriented methods have been used to compute the observability Gramian for nonlinear systems with output measurements, based on partitioning the phase space into slow and fast time regimes using the Lyapunov measure.  
The linear observability Gramian appears here as a special case of the transfer operator based approach. 
The degree of observability of a set can further be related to the residence time. 
These ideas have been extended to determine sensor and actuator placement for the controlled PDE associated with the infinitesimal generator of the transfer operator~\cite{vaidya2012jmaa}.
In particular, the observability and controllability is characterized in terms of the support of the observability and controllability matrices, respectively, and the finite-time and infinite-time Gramians are formulated in terms of the transfer operators. 
The infinite-time controllability is only well defined when the vector field $\bF(\bx)$ is a.e. uniformly stable; thus, an alternative approach to compute the infinite-time controllability set is based on an exponential discounting of the density, which does not rely on the stability property~\cite{sinha2013ecc}.  
The terms \emph{coarse observability} and \emph{coarse controllability} have been introduced to characterize the weaker notion obtained from the finite-dimensional approximation~\cite{sinha2016jmaa}.

The Koopman operator perspective may also be used to assess nonlinear observability and controllability.
Generally, two states are indistinguishable if their future output sequences are identical. Thus, a nonlinear system is considered nonlinearly observable if any pair of states are distinguishable~\cite{HermannandKrener(1977)}.
In~\cite{Surana2016nolcos}, nonlinear observability is evaluated in the proposed Koopman-based linear observer design framework, which provides a linear representation of the underlying system in terms of the Koopman eigenfunctions (see Sec.~\ref{Sec:StateEstimation}). In particular, the nonlinear system is nonlinearly observable if the pair $(\bA,\bC^H)$ is observable, which can be determined via the rank condition of the corresponding observability matrix. 
These ideas have been applied to study pedestrian crowd flow~\cite{benosman2017ifac}, extended further to input-output nonlinear systems~\cite{Surana2016cdc} resulting in bilinear or Lipschitz formulations, and used to compute controllability and reachability~\cite{goswami2017cdc}. 
The observability and controllability Gramians, which are used to examine the degree of observability and controllability, can be computed in the lifted observable space given by a dictionary of functions~\cite{yeung2018acc}.  
In this case, the observability/controllability of the underlying system is related to the observability/controllability of the observables.
The underlying assumption is that the state, input, and output are representable in terms of a few Koopman eigenfunctions, which relies on a suitable choice of observable functions to retain the relevant information.

\subsection{Observer synthesis for state estimation }\label{Sec:StateEstimation}
Observer design is critical for sensor-based estimation of high-dimensional systems.  
The observer is a dynamical system that estimates the current state $\hat{\bx}_k$ from the history of output measurements of the system. 
The Extended Kalman Filter (EKF) is perhaps the most widely used nonlinear observer, benefiting from high performance and a simple formulation based on linearization of the nonlinear state equations. 
However, EKFs are sensitive to noise, do not have convergence guarantees, have limited performance for strongly nonlinear systems, and assume Gaussian noise and disturbances ~\cite{reif2000ieee,chaves2002ejc}. 

Observer design in the Koopman operator framework involves constructing a state estimator for the system
$\bx_{k+1} = \bF (\bx_{k}, \bu_{k}),\quad \by_k= \bH(\bx_{k})$
using a nonlinear transformation of states $\bx$ and outputs $\by$ (and inputs for systems with actuation).  This results in the following estimator dynamical system
\begin{equation}\label{Eqn:KoopmanObserver}
{\bz}_{k+1} = \bA {\bz}_{k},\quad
{\by}_k 		= \bC^y(\bz_{k}),\quad
\bx_k           = \bC^{x}(\bz_k),
\end{equation}
where $\bC^x$ and $\bC^H$ are formed from the Koopman modes and observer Koopman modes.
System~\eqref{Eqn:KoopmanObserver} can be obtained from a spectral decomposition of the Koopman operator, where $\bz$ represent Koopman eigenfunction coordinates, which has been referred to as Koopman observer form~\cite{Surana2016nolcos} (see also Sec.~\ref{Sec:SystemIdentification_KO}).
The resulting system requires that the state and output lie in the span of a subset of eigenfunctions. 
Luenberger or Kalman observers can then be synthesized based on~\eqref{Eqn:KoopmanObserver}, where convergence can be shown under the observability condition based on the Koopman spectral properties, and shown to exhibit increased performance over the extended Kalman filter~\cite{Surana2016nolcos}.
This work has been extended to input-output systems yielding Lipschitz or bilinear observers, resulting from the representation of the unforced system in eigenfunction coordinates~\cite{Surana2016cdc}, and for constrained state estimation based on a moving horizon analogous to model predictive control~\cite{surana2017cdc}.
Building on kernel DMD and the Koopman observer form, a Luenberger Koopman observer has been applied to pedestrian crowd flow to estimate the full state from limited measurements~\cite{benosman2017ifac}.

The probabilistic perspective based on the Perron-Frobenius operator allows one to incorporate uncertainty into observer synthesis. 
Particle filters based on Monte Carlo sampling are commonly used for forecasting and estimation under uncertainty by approximating the PDF with a finite number of samples. 
However, this approach may suffer from the curse of dimensionality~\cite{snyder2008mwr}, particle degeneracy~\cite{arulampalam2002tsp}, and loss of diversity~\cite{ristic2003book}. 
A particular challenge is the accurate estimation of the filter weights.
The Perron-Frobenius and Liouville operator-based approaches have been shown to be advantageous over Monte-Carlo methods~\cite{daum2006non,runolfsson2009acc}.
More recently~\cite{dutta2011jgcd}, a Bayesian formalism involving the Perron-Frobenius operator has been proposed,
which is shown to achieve superior performance over a generic particle filter and bootstrap filter.
This work has been extended to stochastic systems using a Karhunen-Loeve decomposition of the process noise for uncertainty quantification and state estimation~\cite{dutta2012cca,dutta2013acc,dutta2015ieee}.

\subsection{Control design}
Linear models that capture the dynamics of a nonlinear system provide tremendous opportunities for model-based control, making linear techniques readily applicable to nonlinear systems.
In particular, Koopman-based frameworks are amenable to the application of standard linear control theory, in contrast to the probabilistic formulations based on the Perron-Frobenius operator, which require special care.
Koopman models (see Sec.~\ref{Sec:SystemIdentification_KO}) have been increasingly used in combination with optimal control for trajectory steering, e.g. with LQR~\cite{Brunton2016plosone}, SDRE~\cite{Kaiser2017arxiv}, and MPC~\cite{Korda2016arxiv, Peitz2017arxiv,arbabi2018arxiv,KaKuBr2018arxiv}.
Besides theoretical developments, Koopman-based MPC has been increasingly applied in realistic problems, such as power grids~\cite{Korda2018arxiv}, high-dimensional fluid flows~\cite{hanke2018arxiv}, and experimental robotic systems~\cite{abraham2017arxiv,NxR2017video}. 
Although the predictive power of Koopman models may be sensitive to the particular choice of observables and the training data, MPC provides a robust control framework that systematically compensates for model uncertainty by taking into account new measurements (see also Sec.~\ref{Sec:ControlExamples}).
Optimal control formulations have also been considered for switching problems~\cite{sootla2016ifac,sootla2017arxiv,Peitz2017arxiv}.
Based on a global bilinearization, the underlying dynamical system can be stabilized using feedback linearization~\cite{goswami2017cdc}.
A different strategy aims to shape Koopman eigenfunctions directly, referred to as eigenstructure assignment, which has been examined for fluid flow problems~\cite{hemati2017aiaa}.
In contrast to pole placement, which aims to design the eigenvalues of the closed-loop system, eigenstructure assignment aims to additionally modify the structure of the eigenfunctions in a desirable way, which requires additional inputs.
More recently, feedback stabilization for nonlinear systems has been achieved via a control Lyapunov function-based approach based on the bilinear Koopman system in eigenfunction coordinates~\cite{Surana2016cdc}. 
In particular, a convex optimization problem is formulated to search for the control Lyapunov function for the bilinear Koopman system~\cite{huang2019arxiv}.

In many applications it is of interest to shape the distribution of some quantity, e.g. swarms of UAVs or satellites, spin dynamics in nuclear magnetic resonance (NMR) spectroscopy and imaging~\cite{li2009tac}, dispersions in process industries~\cite{wang2001timc}, oil spill in the ocean,  search mission design~\cite{phelps2014automatica}, and beam dynamics~\cite{propoi1994problems,ovsyannikov2006beam}. 
Of particular interest is the shaping of the asymptotic pdf~\cite{forbes2004jpc,guo2005automatica,zhu2012pem}, so that the stationary pdf of the closed-loop system coincides with a desired distribution. 
This relies on the idea that temporal averages can be replaced by spatial analogues under Birkoff's ergodic theorem~\cite{lasota2013book}, after which the mean or variance of a quantity can be controlled by modifying the stationary pdf of the controlled system. 
Optimal control approaches have also been proposed to modify the pdf in the presence of uncertainty in initial conditions or parameters~\cite{ross2016acc,phelps2016sjco}. 
The process of driving one distribution to another one is further intimately related to Monge-Kantorovich optimal transport theory~\cite{monge1781memoire,kantorovich1942dan,elamvazhuthi2016arxiv,halder2014aac}.
In~\cite{halder2014aac}, optimal transport theory has been used to solve an finite-horizon control problem to achieve a desired distribution, where the optimal control vector field is estimated by solving the associated Liouville equation over the finite horizon. Set-oriented and graph-based methods have also been used to study controllability and optimal transport~\cite{elamvazhuthi2016optimal}.

Set-oriented representations of the dynamics (see Sec.~\ref{Sec:SystemIdentification_PFO}) are amenable to optimal control approaches based on dynamic programming~\cite{Bertsekas2005book} for discrete-time, discrete-state, and discrete-action problems. The discretized control problem can be interpreted as a Markov decision problem (MDP).
The optimization problem for the MDP can then be posed as a linear program, which has been demonstrated for communication networks~\cite{alpcan2006cdc}.
An optimal control formulation in terms of an MDP has been formulated for the discretization Liouville equation and demonstrated for navigation in a maze~\cite{kwee2001ab}.
Similar ideas have been used for controlling communication networks~\cite{alpcan2006cdc} and optimizing the asymptotic pdf over the velocity field in high-dimensional fluid flows~\cite{kaiser2017tcfd}.
The inverse Frobenius-Perron problem (IFPP)~\cite{bollt2000ijbc} formulates the control problem for the stationary pdf as a perturbation problem, which aims to find a perturbed transition matrix and the associated perturbed dynamics that give rise to the target density. 
The Markov transition matrix approximating the Perron-Frobenius operator can also be interpreted as a graph, where the nodes correspond the symbols associated with the discrete states and the vertex weights are the transition probabilities facilitating the application of graph-theoretic methods for optimal or minimal-path searches, e.g. in the context of IFPP~\cite{bollt2001ijbc}. 
Stabilization problems have been successfully solved via the Lyapunov measure equation~\cite{vaidya2010ieee,raghunathan2014ieee}, where the problem can be formulated as a linear program, extended for stochastic systems~\cite{das2017acc}, and applied to high-dimensional fluid flows~\cite{vaidya2009cdc}. 
A different perspective is assumed in~\cite{froyland2016siam}, where a convex optimization problem is formulated to determine local perturbations in the context of optimal mixing. 
In particular, the stochastic kernel is designed as a perturbation to the deterministic kernel of the Perron-Frobenius operator, which controls the diffusion and thus the mixing properties of the underlying system. A similar problem has been studied based on the infinitesimal generator~\cite{froyland2016arxiv}, but resulting in a convex quadratic program.
More recent work~\cite{antown2018arxiv} considers the perturbation of the stochastic part of the dynamical system and provides a solution in closed form. 

The success of MPC has also driven its use for controlling probability distributions. 
A nonlinear MPC problem using set-oriented representations~\cite{ohsumi2010ifac} suffers, however, from the curse of dimensionality. 
Thus, the work has been extended using Monte-Carlo and particle filter methods to estimate the pdf, which is referred to as particle MPC~\cite{ohsumi2011ifac}. In particle MPC the Liouville equation enters the MPC formulation as a constraint to control the mean and variance of the state. 
It is worth noting that data-driven PDF prediction using neural nets is increasingly explored due to faster evaluation/prediction of the PDF in contrast to numerically integrating sample points. The loss metric, i.e.\ the error between the true and estimated PDF from NN, can be computed using automatic differentiation of neural nets~\cite{Nakamura2018siam}, building on physics informed neural net via PDE constraints~\cite{raissi2017arxiv_a}, e.g. the Liouville equation.

\subsection{Sensor and actuator placement}
Determining suitable locations of sensors and actuators for data collection and decision-making is a crucial task in any real-world application. Generally, sensor and actuator placement may refer to selecting spatial locations or specific variables, which are used synonymously here. 
Operator-theoretic placement provides a significant opportunity, as they capture global properties, such as meta-stable states, associated with persistent dynamics, and they generalize to nonlinear systems, e.g. for estimating nonlinear observability and controllability. 
Selecting optimal sensor and actuator locations amounts to a combinatorially hard problem. 
Compressed sensing~\cite{candes2006ieee,Donoho2006ieeetit}, sparsity-promoting algorithms employing an $l_1$-penalty term, and sparse sampling techniques such as gappy POD~\cite{everson1995karhunen} have played an increasingly important role in the context of sensor/actuator selection~\cite{Manohar2017csm}. 

In particular, the underlying low-rank structure of the data can be exploited for efficient sensing. 
Thus, compressed sensing and sparsity-promoting algorithms have also been increasingly combined with POD~\cite{Bai2014aiaa} and DMD~\cite{Jovanovic2014pof,Proctor2014epj,Brunton2015jcd,bai2017aiaa}. 
Cluster-based reduced-order models (CROM)~\cite{Kaiser2014jfm}, which provide a coarse approximation of the Perron-Frobenius operator for PDEs, have been combined with the sparse sensor optimization for classification (SSPOC)~\cite{brunton2016siam} algorithm and pivot locations from the QR factorization (see e.g. discrete empirical interpolation methods (DEIM)~\cite{drmac2016siam}) applied to the dominant POD modes to learn optimized sensor locations tailored to a specific model~\cite{kaiser2018jcp}. Probabilistic dynamics are then preserved if the resulting sensing matrix satisfies the restricted isometry property from compressed sensing.
Recently~\cite{manohar2017arxiv}, sensor placement has been explored for multi-scale systems using multi-resolution DMD and DEIM, where sensor locations are tailored to intermittent or transient mode activity.
The control oriented framework DMD with control~\cite{Proctor2016arxiv} has been unified with compressive DMD~\cite{Brunton2015jcd} for compressive system identification~\cite{bai2017aiaa} to identify low-order models from limited input-output data and enable reconstruction of the high-dimensional, interpretable state. This framework, and similarly eDMDc, provide opportunities for optimized placement of sensors and actuators using linear systems techniques. 
A greedy submodular approach was proposed in~\cite{Surana2016slides} based on the Koopman observer form and a mutual information or observability based criterion. 

Sensor placement and actuator placement can also be performed by leveraging the generator PDEs of the Perron-Frobenius and Koopman families.
Controllability and observability Gramians can be generalized for nonlinear systems based on the (controlled) Liouville and adjoint Liouville equation and subsequently used for sensor and actuator placement by maximizing the support (or the $l_2$ norm) of the finite-time Gramians~\cite{vaidya2012jmaa,sinha2013ecc,sinha2016jmaa}.  
The approach utilizes set-oriented methods, where sensors/actuators are placed in certain cells, and the location can be optimized by solving a convex optimization problem~\cite{sinha2013ecc}. 
A greedy heuristic approach based on these ideas using set-oriented methods is proposed in~\cite{fontanini2016be}, which further investigates different criteria such as maximizing the sensing volume (sensor coverage), response time and accuracy (relative measure transported to the sensor in finite time) and incorporating spatial constraints.
The framework has been further extended to incorporate uncertainty~\cite{sharma2018arxiv}. 
Balanced truncation has recently been used for efficient placement of sensors and actuators to simultaneously maximize the controllability and observability Gramians~\cite{Manohar2018arxivB}, which may be promising for Koopman-based Gramians. 

\section{Conclusions \& Future directions}
\label{Sec:Conclusions}

In this chapter, we have explored a number of applications of the Koopman and Perron Frobenius operators for the control of nonlinear systems. 
Although operator-theoretic analysis has a long history in dynamical systems and control, there has been considerable renewed interest in the past decade.  
Much of this recent progress has been driven by the increased availability of data, along with advanced machine learning and optimization algorithms.  
These data-driven approaches are providing approximations to the Koopman and Perron Frobenius operators that are then used for control. 
We have introduced a general control formulation and discussed how this may benefit from a coordinate transformation where nonlinear dynamics become approximately linear.  
Considerations such as stability, controllability and observability, controller and observer design, and sensor/actuator placement were discussed in the context of operator-theoretic approaches, highlighting recent theoretical and applied advances. 
We also discussed a number of important embedding strategies, such as extended DMD, delay coordinates, and eigenfunctions, which may all be used for these control objectives.  
Finally, we provided an example that leverages these various embeddings for model predictive control, showing that MPC is remarkably robust to model uncertainty.  

Despite the tremendous recent progress in operator-based control, there are a number of limitations that must be addressed in the future.  
Obtaining accurate representations of the Koopman operator from data is a key enabler of model-based control. 
However, identifying a coordinate system that is closed under the Koopman operator is notoriously difficult, as the eigenfunctions may have arbitrarily complex representations.  
Fortunately, emerging techniques in deep learning are providing powerful representations for these linearizing transformations~\cite{mardt2017arxiv,otto2017arxiv,takeishi2017anips,lusch2017arxiv,yeung2017arxiv,Li2017chaos,wehmeyer2018arxiv}. 
Another key challenge is the existence of a continuous eigenvalue spectrum in chaotic dynamical systems, which complicates matters considerably.  
Recent methods in delay coordinates~\cite{Susuki2015cdc,Brunton2017natcomm,Arbabi2016arxiv,susuki2017arxiv} and tailored neural network architectures~\cite{lusch2017arxiv} are making headway, although this is far from resolved.  
However, at this point, these advanced embedding techniques have not been adapted for control design, which is an exciting avenue of ongoing work.  

Although much of the focus in data-driven Koopman and Perron Frobenius theory has been on obtaining increasingly accurate approximate embeddings, it may be that improved models have marginal benefit in control design.  
Certainly improved models are useful for prediction and estimation, but advanced control techniques, such as model predictive control, are remarkably robust to lackluster models.  
In many cases, it may be that a simple DMDc model may work nearly as well as a sophisticated model based on the lates embedding techniques.  
However, future operator-based controllers will need to be certified, requiring guarantees on the model fidelity and a rigorous quantification of errors and uncertainties.  
Importantly, the model uncertainty is intimately related to the quality and quantity of the training data, the richness of the dynamics sampled, and the basis used for regression.  
Future efforts will undoubtedly continue to incorporate known or partially known dynamics, structure, and constraints in the data-driven regressions, improving both the models and the uncertainty bounds.  

The current state of the art in operator-theoretic control, with a compelling set of success stories and a long list of future work, make it likely that this field will continue to grow and develop for decades to come.  
It is inevitable that these efforts will increasingly intersect with the fields of machine learning, optimization, and control.  
Advances in machine learning and sparse optimization will continue to drive innovations, both for model discovery and for sensor placement, which are closely related to Koopman and Perron Frobenius.  
Likewise, concrete successes in the control of nonlinear systems will continue to motivate these advanced operator-based approaches.

\begin{acknowledgement}
	EK gratefully acknowledges support by the ``Washington Research Foundation Fund for Innovation in Data-Intensive Discovery" and a Data Science Environments project award from the Gordon and Betty Moore Foundation (Award \#2013-10-29) and the Alfred P. Sloan Foundation (Award \#3835) to the University of Washington eScience Institute, and funding through the Mistletoe Foundation.
	SLB and JNK acknowledge support from the Defense Advanced Research Projects Agency (DARPA contract  PA-18-01-FP-125). SLB acknowledges support from the Army Research Office (W911NF-17-1-0306 and W911NF-17-1-0422). JNK acknowledges support from the Air Force Office of Scientific Research (FA9550-19-1-0011).
\end{acknowledgement}
%
%
%

\bibliographystyle{plain}        
\bibliography{references} 
\end{document}